\documentclass{article}

\usepackage{arxiv}

\usepackage[utf8]{inputenc} 
\usepackage[T1]{fontenc}    
\usepackage{hyperref}       
\usepackage{url}            
\usepackage{booktabs}       
\usepackage{amsfonts}       
\usepackage{nicefrac}       
\usepackage{microtype}      
\usepackage{lipsum}		
\usepackage{graphicx}
\usepackage[numbers]{natbib}
\usepackage{doi}
\usepackage{amsmath}
\usepackage{amssymb}

\graphicspath{{figs/}}

\usepackage{color}
\def\padd#1 {\textcolor{blue}{#1}}

\usepackage{arydshln}
\usepackage{bm}
\makeatletter
\newcounter{letter}
\@whilenum\value{letter}<27\do{
  \expandafter\edef\csname bm\Alph{letter}\endcsname{\noexpand\bm{\Alph{letter}}}
  \expandafter\edef\csname bm\alph{letter}\endcsname{\noexpand\bm{\alph{letter}}}
  \stepcounter{letter}%
}
\newcommand{\bmphi}{\pmb{\phi}}
\newcommand{\bmPhi}{\pmb{\Phi}}
\newcommand{\bmxi}{\pmb{\xi}}
\newcommand{\bmXi}{\pmb{\Xi}}
\newcommand{\bmrho}{\pmb{\rho}}
\newcommand{\bmpsi}{\pmb{\psi}}
\newcommand{\bmPsi}{\pmb{\Psi}}
\newcommand{\bmvarphi}{\pmb{\varphi}}
\newcommand{\bmalpha}{\pmb{\alpha}}

\title{Parallel implementation of a compatible high-order meshless method for the Stokes' equations}

\author{ 
	\hspace{-15mm}\href{https://orcid.org/0000-0001-9246-3122}{\includegraphics[scale=0.06]{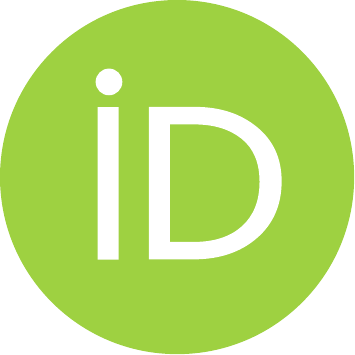}\hspace{1mm}Quang-Thinh Ha} \\
	\hspace{-15mm}Department of Mechanical Engineering\\
	\hspace{-15mm}Boston University\\
	\hspace{-15mm}Boston, MA \\
	\hspace{-15mm}\texttt{qth20@bu.edu} \\
	\And
	\hspace{-5mm}\href{https://orcid.org/0000-0002-2426-4591}{\includegraphics[scale=0.06]{orcid.pdf}\hspace{1mm}Paul A. Kuberry} \\
	\hspace{-5mm}Center for Computing Research\\
	\hspace{-5mm}Sandia National Laboratories\\
	\hspace{-5mm}Albuquerque, NM \\
	\hspace{-5mm}\texttt{pakuber@bu.edu} \\
	\AND
	Nathaniel A. Trask \\
	Center for Computing Research\\
	Sandia National Laboratories\\
	Albuquerque, NM \\
	\texttt{natrask@bu.edu} \\
	\And
    Emily M. Ryan\thanks{Corresponding author: Emily M. Ryan, Department of Mechanical Engineering, Boston University, 110 Cummington Mall, Boston, MA, 02215.  \texttt{ryanem@bu.edu}} \\
	Department of Mechanical Engineering\\
	Boston University\\
	Boston, MA \\
	\texttt{ryanem@bu.edu} \\
}

\hypersetup{
pdftitle={Parallel implementation of a compatible high-order meshless method for the Stokes' equations},
pdfsubject={arXiv},
pdfauthor={Quang-Thinh Ha, Paul A. Kuberry, Nathaniel A. Trask, Emily M. Ryan},
pdfkeywords={meshless, parallel implementation, stokes, moving least squares, staggered scheme, divergence-free, compatible discretization},
}

\begin{document}
\maketitle

\begin{abstract}
	A parallel implementation of a compatible discretization scheme for steady-state Stokes problems is presented in this work. The scheme uses generalized moving least squares to generate differential operators and apply boundary conditions. This meshless scheme allows a high-order convergence for both the velocity and pressure, while also incorporates finite-difference-like sparse discretization. Additionally, the method is inherently scalable: the stencil generation process requires local inversion of matrices amenable to GPU acceleration,  and the divergence-free treatment of velocity replaces the traditional saddle point structure of the global system with elliptic diagonal blocks amenable to algebraic multigrid. The implementation in this work uses a variety of Trilinos packages to exploit this local and global parallelism, and benchmarks demonstrating high-order convergence and weak scalability are provided.
\end{abstract}

\keywords{meshless \and parallel implementation \and stokes \and moving least squares \and staggered scheme \and divergence-free \and compatible discretization}

\section{Introduction}\label{introduction}

Reactive transport in porous media is a multiscale problem relevant to a broad range of applications, including: water filtration systems, carbon sequestration, and pollutant transport \citep{Xiao2012,LANE2018409,LANE20167,ryan2019}. A description of the bulk behavior requires resolution of transport through pores at the submicron level to capture the tortuous transport leading to anomalous diffusion \citep{metzler2014anomalous}. For reactive transport in particular, such a pore-scale resolution is critical to characterizing both the available surface area and exposure time governing chemical reactions \citep{ginn2000distribution,luthy1997sequestration,RYAN201156,RYAN201161}. From a computational perspective however, resolution of the pore scales poses a challenging multiscale problem; as the pathways are topologically complex, the task of generating a high-quality computational mesh is time consuming. Mesh generation in industrial settings is a well-established computational bottleneck in the design-to-analysis process \citep{boggs2005dart}; the added complexity of subsurface pore-scale geometry adds additional complexity beyond comparatively simple manufactured components.

Meshfree discretizations offer an attractive means of automating this process. Previous works have applied methods such as smoothed particle hydrodynamics (SPH) and lattice Boltzmann methods to this end \citep{yang2016intercomparison,tartakovsky2016smoothed}. These techniques typically treat the incompressibility constraint governing hydrodynamics either explicitly, by introducing an artificial compressibility \citep{merkle1987time,morris1997modeling}, or implicitly, by applying pressure projection methods \citep{cummins1999sph}. For the Stokes regime governing flow at such small scales, it is well known that solving the fully coupled Stokes equations is critical to resolving lubrication effects and appropriately treating mass conservation \citep{falk2013stokes,guzman2014conforming}. For meshfree discretizations, solution of the Stokes problem is generally difficult, as the lack of variational principles render a rigorous analysis of inf-sup compatibility impossible without reintroducing a mesh \citep{trask2018compatible}. 

In this work, we develop a parallel implementation of the scheme introduced in \citep{trask2018compatible,hu2019spatially} to provide a scheme appropriate for large scale simulation of fluid flow in Stokesian regimes requiring only easily generated point cloud representations of pores. The scheme uses the generalized moving least squares (GMLS) approximation framework to generate approximations to requisite differential operators \citep{mirzaei2012generalized}. The key idea to avoiding inf-sup instability is to use separate treatment for the velocity and pressure. The velocity is handled via a divergence-free vector polynomial reconstruction while the pressure is treated in a staggered sense \citep{trask2017high}. This allows a decoupling of the velocity/pressure in the Stokes system on the interior of the domain into a block diagonal system, requiring coupling only through boundary terms. 

Despite the challenges in generating stable meshfree discretizations, a few meshfree discretizations of the Stokes problem exist in the literature (see for example \citep{li2020divergence,keim2016high,ortiz2011maximum}). In the current work we have selected \citep{trask2018compatible,hu2019spatially} as a candidate for large scale implementation due to several attractive computational features. Firstly, the GMLS approximation process requires inversion of small dense linear systems at each particle. This local work is ideal for introducing hardware acceleration and allows the generation of high-order stencils with only a modest increase in computational expense. Secondly, the resulting linear system is sparse with positive-definite diagonal blocks, allowing the application of off-the-shelf block preconditioning techniques. To exploit these features, the current work employs the COMpatible PArticle DiscREtization (COMPADRE) toolkit \citep{paul_kuberry_2020_3876465} in Trilinos \citep{heroux2005overview}. COMPADRE provides a means to access: graphics processing unit (GPU) acceleration through the Kokkos library \citep{edwards2014kokkos}, efficient block preconditioners through Teko \citep{cyr2016teko}, fast algebraic multigrid preconditions through MueLU \citep{prokopenko2014muelu}, scalable distribution of particles through Zoltan \citep{boman2012zoltan}, and efficient kd-tree construction of neighbor lists through nanoflann \citep{blanco2014nanoflann}. In concert, these tools provide the necessary ingredients to generate a scalable 3D implementation of the scheme developed in \citep{trask2018compatible}. The focus of this current work is to validate this implementation and demonstrate its strong and weak scaling properties.

The paper is organized as followed. First, a brief review of relevant fundamental of GMLS is presented in Section \ref{gmls}, referring readers interested in further details and applications to \citep{trask2018compatible,hu2019spatially}. Section \ref{applications} establishes how the derivative estimators may be assembled into a collocation scheme for Stokes flow and establish the resulting block matrix structure. Finally, Section \ref{results} provide computational results, using manufactured smooth solutions and analytic solutions to flow past a sphere to establish high-order convergence rates and strong/weak scaling limits of the code.

\section{Generalized Moving Least Squares Approximations}\label{gmls}

GMLS builds a polynomial approximation for a function by minimizing the least-squares errors at specific locations throughout the entire domain. These minimization problems only use points within a small vicinity of the investigating location. The main ingredients of GMLS include:
\begin{itemize}
    \item a set of data points distributed quasi-uniformly throughout the domain,
    \item a set of basis functions, here taken to be polynomials, that represent the approximation and,
    \item the objective function, representing the approximation error, that is minimized in the local least-square problems.
\end{itemize}
By altering the definition of the function space and the local least-square objective function, different variants of the GMLS approximation can be obtained, several of which are used in solving steady state Stokes flow.

\subsection{Basic GMLS approximation}
\label{gmls-basic}

A set of $N$ scattered ``particles'' or ``nodes'' are given at locations $\bmx_i \in \Omega, i=1, \ldots, N$, and the values of the function $u(\bmx) \in C^{\infty} (\Omega)$ at those locations, $u_i = u(\bmx_i)$. The first ingredient for GMLS approximation requires that the distribution of these points is quasi-uniform \citep{Wendland2004}. An approximate reconstruction of the function $u$ from its values $u_i$ is computed by choosing a target point $\bmx_T$, and a polynomial approximation for $u$, valid in the neighborhood of $\bmx_T$, is constructed over the domain.
\begin{equation}
    u_{\text{GMLS}}(\bmx; \bmx_T) = \sum_{i=1}^{Q} p_i(\bmx)^T c_i(\bmx_T) = \bmp(\bmx) \bmc(\bmx_T). \label{eq:polynomial-approximation}
\end{equation}
The approximation in equation (\ref{eq:polynomial-approximation}) is called a local polynomial reproduction of order $m$ if for the family of functions $\{ p_i(\bmx)\} < \infty, (i = 1, 2, \ldots, Q)$, the approximation is exact for all polynomials of order $m$. In equation (\ref{eq:polynomial-approximation}), $\bmp(\bmx) = [p_1(\bmx), p_2(\bmx), \ldots, p_Q(\bmx)]$ is a row vector of linearly independent expansion functions such that $\text{span}(p_1, \ldots, p_Q) = \pi_m(\mathbb{R}^d)$, and $\bmc(\bmx_T) = [c_1, c_2, \ldots, c_Q]^T$ is a column vector of coefficients specific to the target location. This is the second necessity for GMLS: the set of basis functions representing the approximation.

The coefficient column vector $\bmc$ is computed by solving the following optimization problem for a fixed point $\bmx_T$:
\begin{align}
    u_{\text{GMLS}}(\bmx; \bmx_T) &= \min_{\bmc(\bmx_T) \in \mathbb{R}^{Q \times 1}} J_{LS}, \label{eq:optimization-discretized-1}\\
    J_{LS} &= \sum_{j=1}^{N_T} \left[ u_j - u_{\text{GMLS}}(\bmx_j; \bmx_T) \right]^2 W(|| \bmx_T - \bmx_j ||), \label{eq:optimization-discretized-2}\\
    {} &= \sum_{j=1}^{N_T} \left[ u_j - \bmp(\bmx_j)^T \bmc(\bmx_T) \right]^2 W(|| \bmx_T - \bmx_j ||), \label{eq:optimization-discretized-3}\\
    {} &= \sum_{j=1}^{N_T} \left[ u_j - \sum_{i=1}^Q p_i(\bmx_j)c_i(\bmx_T) \right]^2 W(|| \bmx_T - \bmx_j ||),\label{eq:optimization-discretized-4}
\end{align}
where $\bmx_j, j = 1, \ldots, N_T$ are the locations of those particles within a distance $\epsilon$ of $\bmx_T$. Here, $W(r)$ is a radially symmetric kernel with compact support over a ball of radius $\epsilon$. From \citep{Trask2018}, $W(r)$ is chosen to be:
\begin{equation}
    W(r) = \begin{cases}
        \left( 1 - \left( \frac{r}{\epsilon} \right)^4 \right), \quad r < \epsilon, \\
        0, \quad \text{otherwise}.
    \end{cases}
    \label{eq:kernel}
\end{equation}
The error functional $J_{LS}$ in equation (\ref{eq:optimization-discretized-4}) is the sum of the squares of the residuals of all data points within the ``support region'' of the target $\bmx_T$. The three ingredients for GMLS is now fulfilled, and we can proceed with solving for the coefficient vectors $\bmc$.

Equation (\ref{eq:optimization-discretized-4}) can be rewritten as:
\begin{equation}
    J_{LS} = \left( \bmP \bmc(\bmx_T) - \bmu \right)^T \bmW \left( \bmP \bmc(\bmx_T) - \bmu \right), \label{eq:error-functional}
\end{equation}
where:
\begin{align}
    \bmu &= \left[ u_1, u_2, \ldots, u_{N_T} \right]^T = \left[ u(\bmx_1), u(\bmx_2), \ldots, u(\bmx_{N_T}) \right]^T, \label{eq:vector-u} \\
    \bmP &= \begin{bmatrix}
        p_1 (\bmx_1) & p_2(\bmx_1) & \ldots & p_Q(\bmx_1) \\
        p_1 (\bmx_2) & p_2(\bmx_2) & \ldots & p_Q(\bmx_2) \\
        \vdots & \vdots & \ddots & \vdots \\
        p_1 (\bmx_{N_T}) & p_2(\bmx_{N_T}) & \ldots & p_Q(\bmx_{N_T}) \\
    \end{bmatrix} = \begin{bmatrix}
        \bmp(\bmx_1) \\
        \bmp(\bmx_2) \\
        \vdots \\
        \bmp(\bmx_{N_T})
    \end{bmatrix},
    \label{eq:matrix-P} \\
    \bmW &= \text{diag} \left[W(d_1), W(d_2), \ldots, W(d_{N_T}) \right], \quad d_j = ||\bmx_T - \bmx_j||. \label{eq:matrix-W}
\end{align}
Taking the partial derivatives of the cost function $J_{LS}$ in equation (\ref{eq:error-functional}) with respect to the unknown coefficients $\bmc(\bmx_T)$ and equating them to zero gives us:
\begin{align}
    \frac{\partial J_{LS}}{\partial \bmc(\bmx_T)} &= 0, \label{eq:error-functional-derivative}\\
    \bmA(\bmx_T) \bmc(\bmx_T) &= \bmB(\bmx_T) \bmu, \label{eq:local-linear-algebra}
\end{align}
where the matrices $\bmA(\bmx_T) \in \mathbb{R}^{Q \times Q}$ and $\bmB(\bmx_T) \in \mathbb{R}^{Q \times N_T}$ are defined as:
\begin{align}
    \bmA(\bmx_T) &= \sum_{j=1}^{N_T} W(d_j) \bmp(\bmx_j)^T \bmp(\bmx_j) = \bmP^T \bmW \bmP, \label{eq:matrix-A} \\
    \bmB(\bmx_T) &= \left[ W(d_1) \bmp(\bmx_1), W(d_2) \bmp(\bmx_2), \ldots, W(d_{N_T}) \bmp(\bmx_{N_T}) \right] = \bmP^T \bmW. \label{eq:matrix-B}
\end{align}
From equation (\ref{eq:local-linear-algebra}):
\begin{equation}
    \bmc(\bmx_T) = \bmA(\bmx_T)^{-1} \bmB(\bmx_T) \bmu = \bmR_{\bmu \bmu} (\bmx_T) \bmu,
    \label{eq:vector-c}
\end{equation}
where $\bmR_{\bmu \bmu} (\cdot) \in \mathbb{R}^{Q \times N_T}$. The moving least-squares approximation of a scalar field $u$ at location $\bmx_T$ is therefore given by substituting equation (\ref{eq:vector-c}) into equation (\ref{eq:polynomial-approximation}):
\begin{equation}
    u_{\text{GMLS}}(\bmx; \bmx_T) = \bmp(\bmx) \bmc(\bmx_T) = \bmp(\bmx) \bmR_{\bmu \bmu}(\bmx_T) \bmu.
    \label{eq:polynomial-approximation-with-solution}
\end{equation}

For efficiency and simplicity, in our work, the basis functions $\bmp(\bmx)$ are chosen as the Taylor monomials scaled by the kernel support $\epsilon$ for conditioning purposes. In 1D, the monomials are:
\begin{equation}
    p_{\alpha} (x) = \frac{1}{\alpha!} \left( \frac{x - x_T}{\epsilon} \right)^\alpha.
    \label{eq:basis-functions-1d}
\end{equation}
Here $x_T$ is the x-coordinate of the target point in 1D. For 3D, the basis functions are chosen as the tensor products of 1D monomials: 
\begin{equation}
    p_{\alpha} (\bmx) = \left( \frac{1}{{\alpha_x}!} \left( \frac{x - x_T}{\epsilon} \right)^{\alpha_x} \right) \left( \frac{1}{{\alpha_y}!} \left( \frac{y - y_T}{\epsilon} \right)^{\alpha_y} \right) \left( \frac{1}{{\alpha_z}!} \left( \frac{z - z_T}{\epsilon} \right)^{\alpha_z} \right).
    \label{eq:basis-functions-3d}
\end{equation}
Here, $\bmx_T = \{x_T, y_T, z_T\}$ is the coordinate of the target point.

\subsection{GMLS Variant - Constraining the function space of the approximation}
\label{gmls-divfree}

The GMLS can be modified to capture the incompressibility in the solution of velocity field of incompressible Stokes flow. To do so, the set of basis functions can be altered to become the space of divergence-free vector polynomials instead. Let the set of such $m^{\text{th}}$ order polynomials be defined as $\pi^{\text{div}}_m(\mathbb{R}^d)$. For example, in the two-dimensional case with $m=2$:
\begin{equation*}
    \pi^{\text{div}}_2 (\mathbb{R}^2) = \left[
    \begin{pmatrix}
        1 \\ 0
    \end{pmatrix},
    \begin{pmatrix}
        0 \\ 1
    \end{pmatrix},
    \begin{pmatrix}
        y \\ 0
    \end{pmatrix},
    \begin{pmatrix}
        0 \\ x
    \end{pmatrix},
    \begin{pmatrix}
        x \\ -y
    \end{pmatrix},
    \begin{pmatrix}
        y^2 \\ 0
    \end{pmatrix},
    \begin{pmatrix}
        0 \\ x^2
    \end{pmatrix},
    \begin{pmatrix}
        x^2 \\ -2xy
    \end{pmatrix},
    \begin{pmatrix}
        -2xy \\ y^2
    \end{pmatrix}
    \right]
\end{equation*}

The coefficients can then be evaluated similar to the previous section. Let $\bmv(\cdot)$ be a vector field in $\mathbb{R}^d$. A polynomial approximation for the vector function $\bmv$ can then be constructed, again valid in the neighborhood of $\bmx_T$, over the domain:
\begin{equation}
    \bmv_{\text{GMLS}} (\bmx; \bmx_T) = \sum_{i=1}^{Q_d} \bmpsi_i(\bmx) c_i(\bmx_T) = \bmPsi(\bmx) \bmc(\bmx_T), \label{eq:vector-polynomial-approximation}
\end{equation}
where each $\bmpsi_i(\cdot)$ is a divergence-free vector basis function, and $\bmPsi(\bmx) = [\bmpsi_1(\bmx), \bmpsi_2(\bmx), \ldots, \bmpsi_{Q_d}(\bmx)]$ is a $d-$dimensional vector of $Q_d$ columns and $d$ rows such that $\text{span}(\bmpsi_1, \ldots, \bmpsi_{Q_d}) = \pi^{\text{div}}_m(\mathbb{R}^d)$. The coefficient column vector $\bmc(\cdot) \in \mathbb{R}^{Q_d \times 1}$ is similarly computed by solving an optimization problem for a fixed point $\bmx_T$:
\begin{align}
    \bmv_{\text{GMLS}} (\bmx; \bmx_T) &= \min_{\bmc(\bmx_T) \in \mathbb{R}^{Q_d \times 1}} J_{LS-V}, \label{eq:vector-optimization-discretized-1}\\
    J_{LS-V} &= \sum_{j=1}^{N_T} \left[ \bmv_j - \bmv_{\text{GMLS}}(\bmx_j; \bmx_T) \right] \cdot \left[ \bmv_j - \bmv_{\text{GMLS}}(\bmx_j; \bmx_T) \right] W(d_j), \label{eq:vector-optimization-discretized-2}\\
    {} &= \sum_{j=1}^{N_T} \left[ \bmv_j - \bmPsi(\bmx_j) \bmc(\bmx_T) \right]\cdot \left[ \bmv_j - \bmPsi(\bmx_j) \bmc(\bmx_T) \right] W(d_j), \label{eq:vector-optimization-discretized-3}\\
    {} &= \sum_{j=1}^{N_T} \left[ \bmv_j - \sum_{i=1}^{Q_d} \bmpsi_i(\bmx_j)c_i(\bmx_T) \right] \cdot \left[ \bmv_j - \sum_{i=1}^{Q_d} \bmpsi_i(\bmx_j)c_i(\bmx_T) \right]W(d_j).\label{eq:vector-optimization-discretized-4}
\end{align}
Taking the derivative of equation (\ref{eq:vector-optimization-discretized-4}) with respect to the coefficient vectors:
\begin{equation}
    \frac{\partial J_{LS-V}}{\partial c(\bmx_T)} = -2 \sum_{j = 1}^{N_T} W(d_j) \bmpsi_m (\bmx_j) \cdot \left[ \bmx_j - \sum_{i=1}^{Q_d} \bmpsi_i(\bmx_j) c_i (\bmx_T) \right], \label{eq:derivative-vector-objective-function}
\end{equation}
where the index $m$ in equation (\ref{eq:derivative-vector-objective-function}) loops through each of the basis functions in $\bmPsi(\cdot)$, i.e. $m = 1, \ldots, Q_d$. The solution to the minimization problem in equation (\ref{eq:vector-optimization-discretized-1}) is obtained by setting equation (\ref{eq:derivative-vector-objective-function}) to zero. A straightforward manipulation leads to:
\begin{equation}
    \sum_{j-1}^{N_T} W(d_j) \bmpsi_m (\bmx_j) \cdot \sum_{i=1}^{Q_d} \bmpsi_i (\bmx_j) c_i (\bmx_T) = \sum_{j=1}^{N_T} W(d_j) \bmpsi_m (\bmx_j) \cdot \bmv_j. \label{eq:vector-linear-algebra}
\end{equation}
There are $Q_d$ equations coming from (\ref{eq:vector-linear-algebra}) since $m$ is the only free index. Thus, a linear algebra system of $Q_d$ equations can be obtained to solve for $Q_d$ unknowns of the coefficient vector $\bmc(\bmx_T)$. Finally, the following approximation for the vector function $\bmv_{\text{GMLS}}(\cdot)$, and subsequently for the polynomial approximation of $\bmv_{\text{GMLS}}(\cdot)$, can be obtained:
\begin{align}
    \bmV &= [ \bmv_1, \bmv_2, \ldots, \bmv_{N_T}]^T =
        \begin{bmatrix}
            \vline & \vline & {} & \vline \\
            \bmv (\bmx_1) & \bmv (\bmx_2) & \cdots & \bmv(\bmx_{N_T}) \\
            \vline & \vline & {} & \vline
        \end{bmatrix}^T \label{eq:vector-V} \\
    \bmc(\bmx_T) &= \bmR_{\bmv \bmv} (\bmx_T) \bmV \label{eq:vector-c-2} \\
    \bmv_{\text{GMLS}} (\bmx; \bmx_T) &= \bmPsi(\bmx)^T \bmc(\bmx_T) = \bmPsi(\bmx) \bmR_{\bmv \bmv} (\bmx_T) \bmV \label{eq:vector-solved-approximation}
\end{align}
For the sake of clarity, the dimensions of these matrices and vectors are $\bmR_{\bmv \bmv} (\cdot) \in \mathbb{R}^{Q_d \times N_T}$ and $\bmV (\cdot) \in \mathbb{R}^{N_T \times d}$.

\subsection{GMLS Variant - Incorporate constraints on the local least-square problems}
\label{gmls-neumann}

Another possible variant for GMLS is changing the objective function in the local least-squares problems. By introducing additional constraints, conditions can be imposed on the approximation to, for example, enforce Neumann conditions on boundary particles of a bounded domain $\Omega$. Specifically, this formulation focuses on the following case of:
\begin{equation}
    \bmn \cdot \nabla \bmu = h, \text{ if } \bmx \in \Gamma = \partial \Omega, \label{eq:neumann-l2}
\end{equation}
in which $\bmn$ is the outward normal vector of the domain $\Omega$. The condition (\ref{eq:neumann-l2}) is then enforced through the choice of GMLS expansion coefficients rather than through collocation. Therefore, in order to enforce the constraint in equation (\ref{eq:neumann-l2}), the optimization problem in equation (\ref{eq:optimization-discretized-1})-(\ref{eq:optimization-discretized-4}) is modified as follows:
\begin{align}
    u_{\text{GMLS}} (\bmx; \bmx_T) &= \min_{\bmc(\bmx_T) \in \mathbb{R}^d} J_{LS}, \label{eq:neumann-optimization-1} \\
    J_{LS} &= \sum_{j=1}^{N_T} \left[ u_j - u_{\text{GMLS}} (\bmx_j; \bmx_T) \right]^2 W(|| \bmx_T - \bmx_j ||), \label{eq:neumann-optimization-2} \\
    \text{s.t. } & \bmn \cdot \nabla u_{\text{GMLS}} (\bmx_T; \bmx_T) = h(\bmx_T).\label{eq:neumann-optimization-3}
\end{align}
This applies only to target points on the Neumann boundary. The constrained optimization problem in equation (\ref{eq:neumann-optimization-1})-(\ref{eq:neumann-optimization-3}) can be rewritten to include a Lagrange multiplier. Using the polynomial approximation in equation (\ref{eq:polynomial-approximation}) and following a similar procedure yields:
\begin{align}
    u_{\text{GMLS}} (\bmx; \bmx_T) &= \min_{\bmc(\bmx_T) \in \mathbb{R}^d} J_{LS-N}, \label{eq:neumann-optimization-discretized-1} \\
    J_{LS-N} &= \sum_{j=1}^{N_T} \left[ u_j - u_{\text{GMLS}} (\bmx_j; \bmx_T) \right]^2 W(|| \bmx_T - \bmx_j ||) \nonumber \\
    {} & \qquad + \lambda(\bmx_T) \left( \bmn(\bmx_T) \cdot \nabla u_{\text{GMLS}} (\bmx; \bmx_T) - h(\bmx_T) \right), \label{eq:neumann-optimization-discretized-2} \\
    {} &= \sum_{j=1}^{N_T} \left[ u_j - \bmp(\bmx_j) \bmc(\bmx_T) \right]^2 W(|| \bmx_T - \bmx_j ||) \nonumber \\
    {} & \qquad + \lambda(\bmx_T) \left[ \left( \left(\bmn(\bmx_T) \cdot \nabla \bmp(\bmx)|_{\bmx=\bmx_T} \right)^T \bmc(\bmx_T) \right) - h(\bmx_T) \right], \label{eq:neumann-optimization-discretized-3}\\
    {} &= \sum_{j=1}^{N_T} \left[ u_j - \sum_{i=1}^Q p_i(\bmx_j)c_i(\bmx_T) \right]^2 W(|| \bmx_T - \bmx_j ||) \nonumber \\
    {} & \qquad + \lambda(\bmx_T) \left[ \sum_{k=1}^3 n_k (\bmx_T) \left( \sum_{i=1}^Q \left( \nabla p_i(\bmx)|_{\bmx = \bmx_T} \right) c_i (\bmx_T)  \right)_k - h(\bmx_T) \right]. \label{eq:neumann-optimization-discretized-4}
\end{align}
It is noted that in equation (\ref{eq:neumann-optimization-discretized-4}), the subscript $k$ denotes the spatial component of a vector field. Additionally, when taking the derivative of the polynomial approximation $u_{\text{GMLS}}$, the spatial dependence of the coefficient vector $\bmc$ is neglected. Similar to before, the functional $J_{LS-N}$ can be consolidated into matrix-vector multiplication form:
\begin{align}
    J_{LS-N} &= \left( \bmP \bmc(\bmx_T) - \bmu \right)^T \bmW \left( \bmP \bmc(\bmx_T) - \bmu \right) \nonumber \\
    {} & \qquad + \lambda(\bmx_T) \left[ \left( \bmn(\bmx_T) \cdot \left( \nabla \bmp(\bmx)|_{\bmx=\bmx_T} \right) \bmc(\bmx_T) \right) - h(\bmx_T) \right].
    \label{eq:neumann-error-functional}
\end{align}
Here the gradient of the expansion functions takes the form of:
\begin{equation}
    \nabla \bmp(\bmx) = \begin{bmatrix}
        \frac{\partial p_1 (\bmx)}{\partial x} & \frac{\partial p_2 (\bmx)}{\partial x} & \ldots & \frac{\partial p_Q (\bmx)}{\partial x}\\
        \frac{\partial p_1 (\bmx)}{\partial y} & \frac{\partial p_2 (\bmx)}{\partial y} & \ldots & \frac{\partial p_Q (\bmx)}{\partial y}\\
        \frac{\partial p_1 (\bmx)}{\partial z} & \frac{\partial p_2 (\bmx)}{\partial z} & \ldots & \frac{\partial p_Q (\bmx)}{\partial z}\\
    \end{bmatrix}
\end{equation}

Taking the derivative of $J_{LS-N}$ in equation (\ref{eq:neumann-error-functional}) with respect to the unknown coefficient vector $\bmc(\bmx_T)$ and the Lagrange multiplier $\lambda(\bmx_T)$ gives the following saddle-point system:
\begin{align}
    \frac{\partial J_{LS-N}}{\partial \bmc(\bmx_T)} &= 0, \label{eq:neumann-error-functional-derivative-c}\\
    \frac{\partial J_{LS-N}}{\partial \lambda(\bmx_T)} &= 0, \label{eq:neumann-error-functional-derivative-lambda}\\
    \left[
        \begin{array}{c;{2pt/2pt}c}
            \bmP^T \bmW \bmP & \bmF(\bmx_T)  \\ \hdashline[2pt/2pt]
            \bmF(\bmx_T)^T & 0
        \end{array}
    \right]
    \left[
        \begin{array}{c}
            \bmc(\bmx_T) \\ \hdashline[2pt/2pt]
            \lambda
         \end{array}
    \right]
    &= \left[
        \begin{array}{c}
            \bmP^T \bmW \bmu \\ \hdashline[2pt/2pt]
            h(\bmx_T)
        \end{array}
    \right]
    = \left[
        \begin{array}{c;{2pt/2pt}c}
            \bmP^T \bmW & \mathbf{0}  \\ \hdashline[2pt/2pt]
            \mathbf{0} & 1
        \end{array}
    \right]
    \left[
        \begin{array}{c}
            \bmu \\ \hdashline[2pt/2pt]
            h(\bmx_T)
         \end{array}
    \right], \label{eq:neumann-saddle-point-system}
\end{align}
where $\bmF(\bmx_T) = \nabla \bmp(\bmx_T)^T \bmn(\bmx_T)$, and $\bmu$, $\bmP$ and $\bmW$ are defined in equations (\ref{eq:vector-u}), (\ref{eq:matrix-P}) and (\ref{eq:matrix-W}), respectively. By solving the linear algebra problem in equation (\ref{eq:neumann-saddle-point-system}), a new expression for the coefficient vector $\bmc(\bmx_T)$ can be obtained:
\begin{align}
    \left[
        \begin{array}{c;{2pt/2pt}c}
            \bmR_{\bmu \bmu}(\bmx_T) & R_{\bmu h}(\bmx_T)  \\ \hdashline[2pt/2pt]
            \bmR_{h \bmu}(\bmx_T) & R_{hh}(\bmx_T)
        \end{array}
    \right] &=
    \left[
        \begin{array}{c;{2pt/2pt}c}
            \bmP^T \bmW \bmP & \bmF(\bmx_T)  \\ \hdashline[2pt/2pt]
            \bmF(\bmx_T)^T & 0
        \end{array}
    \right]^{-1}
    \left[
        \begin{array}{c;{2pt/2pt}c}
            \bmP^T \bmW & \mathbf{0}  \\ \hdashline[2pt/2pt]
            \mathbf{0} & 1
        \end{array}
    \right], \label{eq:neumann-kkt-linear-algebra} \\
    \left[
        \begin{array}{c}
            \bmc(\bmx_T) \\ \hdashline[2pt/2pt]
            \lambda
         \end{array}
    \right] &=
    \left[
        \begin{array}{c;{2pt/2pt}c}
            \bmR_{\bmu \bmu} (\bmx_T) & \bmR_{\bmu h} (\bmx_T) \\ \hdashline[2pt/2pt]
            \bmR_{h \bmu} (\bmx_T)& R_{hh}(\bmx_T)
        \end{array}
    \right]
    \left[
        \begin{array}{c}
            \bmu \\ \hdashline[2pt/2pt]
            h(\bmx_T)
         \end{array}
    \right], \label{eq:neumann-solution-coefficient-lambda}
\end{align}
or, with $\bmR_{\bmu \bmu} (\cdot) \in \mathbb{R}^{Q \times N_T}$ and $\bmR_{\bmu h} (\cdot) \in \mathbb{R}^{Q \times 1}$:
\begin{equation}
    \bmc(\bmx_T) = \bmR_{\bmu \bmu}(\bmx_T) \bmu + \bmR_{\bmu h}(\bmx_T) h(\bmx_T). \label{eq:neumann-solution-coefficient}
\end{equation}

\subsection{GMLS Variant - Derivative approximations by modified objective function}
\label{gmls-staggered}

The last GMLS variant that will be covered in this study is the so-called staggered GMLS scheme. For this case, both the set of basis functions and the local least-squares objective functions are changed. For more details with regards to the implementation and analysis of the staggered scheme, please refer to \citep{Trask2017}.

The staggered scheme approximates a scalar field by approximating its derivative at the midpoint between nodes with a vector function. Analogous to the previous approximations, let $\pi_m^{\bmv} (\mathbb{R}^d)$ be the set of vector polynomials of order at most $m$ and a gradient vector field $\bmq = \nabla p$. The polynomial approximation of $\bmq \approx \bmq_h$ is computed as:
\begin{equation}
    \bmq_h(\bmx; \bmx_T) = \sum_{i=1}^{Q_v} \bmphi_i (\bmx) c_i (\bmx_T) = \bmPhi(\bmx) \bmc (\bmx_T), \label{eq:staggered-polynomial-approximation}
\end{equation}
where $\bmPhi(\bmx) = [\bmphi_1(\bmx), \bmphi_2(\bmx), \ldots, \bmphi_{Q_v}]$ is a collection where each of the column is a vector basis function $\bmphi_i(\bmx)$. It is also required that $\text{span}(\bmpsi_1, \ldots, \bmpsi_{Q_v}) = \pi^{\bmv}_m(\mathbb{R}^d)$. The coefficient vector $\bmc$ is now the solution of the following minimization problem:
\begin{align}
    \bmq_h (\bmx; \bmx_T) &= \min_{\bmc(\bmx_T) \in \mathbb{R}^{Q_v \times 1}} J_{LS-S}, \label{eq:staggered-optimization-discretized-1}\\
    J_{LS-S} &= \sum_{j=1}^{N_T} \left[ \int_{\bmx_T}^{\bmx_j} \bmq_h (\bmx; \bmx_T) \cdot d\bmx - \int_{\bmx_T}^{\bmx_j} \bmq (\bmx) \cdot d\bmx \right]^2  W(d_j). \label{eq:staggered-optimization-discretized-2}
\end{align}
Using the fact that $\bmq = \nabla p$, the term $\int_{\bmx_T}^{\bmx_j} \bmq (\bmx) \cdot d\bmx$ is equal to $p_j - p(\bmx_T)$ via the fundamental theorem of calculus. By substituting equation (\ref{eq:staggered-polynomial-approximation}) into the term $\int_{\bmx_T}^{\bmx_j} \bmq_h (\bmx; \bmx_T) \cdot d\bmx$:
\begin{equation}
    \int_{\bmx_T}^{\bmx_j} \bmq_h (\bmx; \bmx_T) \cdot d\bmx = \int_{\bmx_T}^{\bmx_j} \sum_{i=1}^{Q_v} \bmphi_i(\bmx) c_i (\bmx_T) \cdot d\bmx = \sum_{i=1}^{Q_v} c_i (\bmx_T) \int_{\bme_{ij}} \bmphi_i(\bmx) \cdot d\bmx. \label{eq:staggered-integral-term}
\end{equation}
Using equation (\ref{eq:staggered-integral-term}), the objective function in equation (\ref{eq:staggered-optimization-discretized-3}) is now equivalent to:
\begin{equation}
    J_{LS-S} = \sum_{j=1}^{N_T} \left[\left( \sum_{i=1}^{Q_v} c_i (\bmx_T) \int_{\bme_{ij}} \bmphi_i(\bmx) \cdot d\bmx \right) - (p_j - p(\bmx_T)) \right]^2  W(d_j). \label{eq:staggered-optimization-discretized-3}
\end{equation}
The following vectors are defined:
\begin{align}
    \bmxi(\bmx_j)  &= [\xi_1(\bmx_j), \xi_2(\bmx_j), \ldots, \xi_{Q_v}(\bmx_j)], \quad \xi_i(\bmx_j)  = \int_{\bmx_T}^{\bmx_j} \bmphi_i (\bmx)\cdot d \bmx, \label{eq:vector-xi} \\
    \bmXi &= \begin{bmatrix}
        \xi_1(\bmx_1) & \xi_2(\bmx_1) & \ldots & \xi_{Q_v}(\bmx_1) \\
        \xi_1(\bmx_2) & \xi_2(\bmx_2) & \ldots & \xi_{Q_v}(\bmx_2) \\
        \vdots & \vdots & \ddots & \vdots \\
        \xi_1(\bmx_{N_T}) & \xi_2(\bmx_{N_T}) & \ldots & \xi_{Q_v}(\bmx_{N_T})
    \end{bmatrix} = \begin{bmatrix}
        \bmxi(\bmx_1) \\
        \bmxi(\bmx_2) \\
        \vdots \\
        \bmxi(\bmx_{N_T})
    \end{bmatrix}, \label{eq:matrix-Xi} \\
    \bmp &= [p_1, p_2, \ldots, p_{N_T}]^T, \label{eq:vector-p} \\
    \bmrho &= [ \left(p_1 - p(\bmx_{T}) \right), \left( p_2 - p(\bmx_{T}) \right), \ldots, \left( p_{N_T} - p(\bmx_{T}) \right)]^T = \bmD \bmp, \label{eq:vector-rho}
\end{align}
where $\bmD$ is similar to a finite difference stencil. Equation (\ref{eq:staggered-optimization-discretized-3}) can be rewritten into the following form:
\begin{equation}
    J_{LS-S} = \left( \bmXi \bmc(\bmx_T) - \bmrho \right)^T \bmW \left( \bmXi \bmc(\bmx_T) - \bmrho \right). \label{eq:staggered-matrix-vector-multiplication}
\end{equation}
The unknown coefficients can now be found by following now familiar steps:
\begin{align}
    \frac{\partial J_{LS-S}}{\partial \bmc(\bmx_T)} &= 0, \label{eq:staggerd-error-functional-derivative} \\
    \bmA_s(\bmx_T) \bmc(\bmx_T) &= \bmB_s(\bmx_T) \bmrho, \label{eq:staggered-local-linear-algebra} \\
    \bmA_s(\bmx_T) &= \sum_{j=1}^{N_T} W(d_j) \bmxi(\bmx_j)^T \bmxi(\bmx_j) = \bmXi^T \bmW \bmXi, \label{eq:staggerd-matrix-A} \\
    \bmB_s(\bmx_T) &= \left[ W(d_1) \bmxi(\bmx_1), W(d_2) \bmxi(\bmx_2), \ldots, W(d_{N_T}) \bmxi(\bmx_{N_T}) \right] = \bmXi^T \bmW, \label{eq:staggered-matrix-B} \\
    \bmc(\bmx_T) &= \bmA_s(\bmx_T)^{-1} \bmB_s(\bmx_T) \bmrho = \bmA_s(\bmx_T)^{-1} \bmB_s(\bmx_T) \bmD \bmp. \label{eq:staggered-vector-c}
\end{align}
Hence, the vector function $\bmq$, and thus $\nabla p$, can be approximated by substituting equation (\ref{eq:staggered-vector-c}) into equation (\ref{eq:staggered-polynomial-approximation}):
\begin{align}
    \nabla p = \bmq &\approx \bmq_h(\bmx; \bmx_T) \nonumber \\
    {} &= \bmPhi(\bmx) \bmc (\bmx_T) \nonumber \\
    {} &= \bmPhi(\bmx) \bmA_s(\bmx_T)^{-1} \bmB_s(\bmx_T) \bmD \bmp \nonumber \\
    \nabla p &\approx \bmPhi(\bmx) \bmR_{\bmq \bmq} (\bmx_T) \bmp. \label{eq:staggered-polynomial-approximation-solution}
\end{align}
If desired, the Laplacian of $p$ can now be approximated by taking divergence on both sides of equation (\ref{eq:staggered-polynomial-approximation-solution}):
\begin{align}
    \nabla^2 p = \nabla \cdot (\nabla p) = \nabla \cdot \bmq &\approx \nabla \cdot \bmq_h(\bmx; \bmx_T) \nonumber \\
    {} &= \nabla \cdot \left( \bmPhi(\bmx) \bmR_{\bmq \bmq} (\bmx_T) \bmp \right) \nonumber \\
    {} &= \left( \nabla \cdot \bmPhi(\bmx) \right) \bmR_{\bmq \bmq} (\bmx_T) \bmp \label{eq:staggered-laplacian-approximation}
\end{align}

\subsection{GMLS approximation of operators}

Once the approximation of the field $u(\bmx)$ at the target location $\bmx_T$ is obtained, the operators on $\bmu(\cdot)$ at said location can also be approximated. For example, assuming $\Omega \subset \mathbb{R}^n$ and $u(\bmx) = u(x_1, \ldots, x_n)$, using multi-index notation of:
\begin{align*}
    |\beta| &= \beta_1 + \beta_2 + \cdots + \beta_n \\
    D^{\beta} u(\bmx) &= \frac{ \partial^{|\beta|} }{\partial x_1^{\beta_1} \cdots \partial x_n^{\beta_n}} u(x_1, \ldots, x_n)
\end{align*}
the approximation of the derivatives of $u(\bmx)$ at $\bmx_T$ is obtained by differentiating equation (\ref{eq:polynomial-approximation}) with respect to $\bmx$:
\begin{equation}
    D^{\beta} u(\bmx) \approx D^{\beta} u_{\text{GMLS}}(\bmx; \bmx_T) = D^{\beta} \bmH(\bmx) \bmu,
    \label{eq:polynomial-approximation-derivative-concept}
\end{equation}
Noticeably, the spatial dependence of $\bmc$ in equation (\ref{eq:polynomial-approximation-derivative}) is neglected, and the derivatives are applied directly to the polynomial basis $\bmp(\bmx)$. Thus:
\begin{equation}
    D^{\beta} u(\bmx) \approx D^{\beta} \left( \bmH(\bmx) \bmu \right) = \left( D^{\beta} \bmp(\bmx) \right) \bmR_{\bmu \bmu} (\bmx_T) \bmu.
    \label{eq:polynomial-approximation-derivative}
\end{equation}
Previous work shows that this approximation is bounded and converges to the exact derivative with a rate of $ m + 1 - \left| \beta \right|$ \citep{Armentano2001}.

\section{Applications - Stokes flow}\label{applications}

\subsection{Formulation}

In this work, the steady Stokes problem is considered:
\begin{alignat}{2}
    -\nu \nabla^2 \bmv + \nabla \varphi &= \bmf, \quad &&\bmx \in \Omega, \label{eq:stokes-momentum-conservation} \\
    \nabla \cdot \bmv &= 0, \quad &&\bmx \in \Omega \label{eq:stokes-mass-conservation} \\
    \bmv &= \bmw, \quad &&\bmx \in \Gamma, \label{eq:stokes-velocity-boundary-condition}
\end{alignat}
where the domain $\Omega \subset \mathbb{R}^d$ has a piecewise continuous boundary $\Gamma = \partial \Omega$, $\bmv(\bmx)$ and $\varphi(\bmx)$ are velocity and pressure, $\nu$ is the kinematic viscosity, and $\bmf(\bmx)$ and $\bmw(\bmx)$ are given data and Dirichlet conditions.

Using the vector identity $\nabla ^2 \bmv = - \nabla \times \nabla \times \bmv + \nabla (\nabla \cdot \bmv)$ along with equation (\ref{eq:stokes-mass-conservation}), equation (\ref{eq:stokes-momentum-conservation}) can be rewritten as:
\begin{equation}
    \nu \nabla \times \nabla \times \bmv + \nabla \varphi = \bmf, \quad \bmx \in \Omega
    \label{eq:stokes-velocity-divergence-free-vector-alone}
\end{equation}
Taking the divergence on both sides of equation (\ref{eq:stokes-velocity-divergence-free-vector-alone}), the following equation for $\varphi$ is obtained:
\begin{equation}
    \nabla^2 \varphi = \nabla \cdot \bmf, \quad \bmx \in \Omega. \label{eq:stokes-pressure-laplacian-alone}
\end{equation}
Now a boundary condition for the pressure field is required. By examining equation (\ref{eq:stokes-velocity-divergence-free-vector-alone}) along the boundary $\Gamma$ of the domain, a dot product of both sides with the unit normal vector $\hat{\bmn}$, which is assumed to be provided, can be performed along the boundary:
\begin{equation}
    \nu \hat{\bmn} \cdot \left( \nabla \times \nabla \times \bmv \right) + \hat{\bmn} \cdot \nabla \varphi = \hat{\bmn} \cdot \bmf, \quad \bmx \in \Gamma \label{eq:stokes-pressure-neumann-alone}.
\end{equation}
Therefore, by choosing the velocity from an appropriate space of divergence-free vector fields so that equation (\ref{eq:stokes-mass-conservation}) holds, the Stokes problem in equations (\ref{eq:stokes-momentum-conservation}) and (\ref{eq:stokes-velocity-boundary-condition}) is equivalent to the following coupled system of equations:
\begin{alignat}{2}
    \nu \nabla \times \nabla \times \bmv + \nabla \varphi &= \bmf, \quad &&\bmx \in \Omega, \label{eq:stokes-velocity-divergence-free-vector} \\
    \bmv &= \bmw, \quad &&\bmx \in \Gamma, \label{eq:stokes-velocity-dirichlet} \\
    \nabla^2 \varphi &= \nabla \cdot \bmf, \quad &&\bmx \in \Omega, \label{eq:stokes-pressure-laplacian} \\
    \nu \hat{\bmn} \cdot \left( \nabla \times \nabla \times \bmv \right) + \hat{\bmn} \cdot \nabla \varphi &= \hat{\bmn} \cdot \bmf, \quad &&\bmx \in \Gamma \label{eq:stokes-pressure-neumann}.
\end{alignat}
Here, the velocity field has a Dirichlet boundary condition applied along the boundary $\Gamma$ in equation (\ref{eq:stokes-velocity-dirichlet}), while the pressure field has a pure Neumann boundary condition from equation (\ref{eq:stokes-pressure-neumann}).

\subsection{Solution of steady state Stokes flow via collocation}
The method used here is described in further details in \citep{Trask2018}. The Stokes collocation formula is solved by enforcing the governing equations (\ref{eq:stokes-velocity-divergence-free-vector}) and (\ref{eq:stokes-pressure-laplacian}) at a certain number of locations $N_\Omega$ inside the domain $\Omega$, while the boundary conditions of equation (\ref{eq:stokes-velocity-dirichlet}) and (\ref{eq:stokes-pressure-neumann}) are being enforced at $N_\Gamma$ points along the boundary $\Gamma = \partial \Omega$. Let $I_\Omega$ and $I_\Gamma$ be the sets of unique internal points and boundary points, respectively. The size of these sets, therefore, are $|I_\Omega| = N_\Omega$ and $|I_\Gamma| = N_\Gamma$, in which the union of these sets is the entire domain of $N$ particles (i.e. $|I_\Omega \cup I_\Gamma| = N_\Omega + N_\Gamma$ and $I_\Omega \cap I_\Gamma = \varnothing$). The governing equations (\ref{eq:stokes-velocity-divergence-free-vector})-(\ref{eq:stokes-pressure-neumann}) can then be specified as followed:
\begin{alignat}{2}
    \nu \nabla \times \nabla \times \bmv(\bmx_I) + \nabla \varphi (\bmx_I) &= \bmf(\bmx_I), \quad && I \in I_\Omega, \label{eq:stokes-velocity-target} \\
    \bmv(\bmx_I) &= \bmw(\bmx_I), \quad && I \in I_\Gamma, \label{eq:stokes-velocity-dirichlet-target} \\
    \nabla^2 \varphi (\bmx_I) &= \nabla \cdot \bmf(\bmx_I), \quad && I \in I_\Omega, \label{eq:stokes-pressure-laplacian-target} \\
    \nu \hat{\bmn} \cdot \left( \nabla \times \nabla \times \bmv(\bmx_I) \right) + \hat{\bmn} \cdot \nabla \varphi(\bmx_I) &= \hat{\bmn} \cdot \bmf(\bmx_I), \quad && I \in I_\Gamma \label{eq:stokes-pressure-neumann-target}.
\end{alignat}

From equation (\ref{eq:vector-solved-approximation}), the velocity field $\bmv$ can be approximated with GMLS using a divergence-free basis. The approximation of the resulting vector function when an operator is applied on it can then be evaluated. For example, the approximate curl on the vector field, at location $\bmx_I$, can be obtained as:
\begin{equation}
    \left( \nabla \times \nabla \times \bmv \right) (\bmx) \approx \nabla \times \nabla \times \bmv_{\text{GMLS}} (\bmx; \bmx_I) = \nabla \times \nabla \times \left( \bmPsi(\bmx) \bmc(\bmx_I) \right), \label{eq:curl-vector-approximation-concept}
\end{equation}
in which, again, the spatial dependence in $\bmc(\cdot)$ is neglected. Equation (\ref{eq:curl-vector-approximation-concept}) can be rewritten with a finite-difference like stencil:
\begin{equation}
    \left( \nabla \times \nabla \times \bmv \right) (\bmx_I) \approx \sum_{j \in \text{supp}(W_{Ij})} \bmalpha^3_{Ij} \bmv_j, \quad \text{ for } I \in \Omega. \label{eq:curl-vector-stencil}
\end{equation}
where $j \in \text{supp}(W_{Ij})$ identifies neighboring particle $j$ of particle $I$.

For particles inside the domain $\Omega$, the Laplacian of the pressure field can be approximated using the staggered scheme described in section \ref{gmls-staggered}. Similarly, the approximation from equation (\ref{eq:staggered-laplacian-approximation}) can be rewritten into another finite-difference stencil:
\begin{equation}
    \nabla^2 \varphi (\bmx_I) \approx \sum_{j \in \text{supp}(W_{Ij})} \bmalpha^2_{Ij} \varphi_j, \quad \text{ for } I \in \Omega. \label{eq:laplacian-pressure-stencil} \\
\end{equation}
In order to incorporate the Neumann boundary conditions for the pressure field, the derivation in section \ref{gmls-neumann} is used. In fact, if equation (\ref{eq:stokes-pressure-neumann}) is rewritten and combined with equation (\ref{eq:stokes-pressure-laplacian}):
\begin{alignat}{2}
    \nabla^2 \varphi &= \nabla \cdot \bmf, \quad &&\bmx \in \Omega, \\
    \hat{\bmn} \cdot \nabla \varphi &= \hat{\bmn} \cdot \bmf - \nu \hat{\bmn} \cdot \left( \nabla \times \nabla \times \bmv \right), \quad &&\bmx \in \Gamma \label{eq:stokes-pressure-neumann-1},
\end{alignat}
then equation (\ref{eq:stokes-pressure-neumann-1}) is equivalent to equation (\ref{eq:neumann-l2}). Let $s_I = \nabla \cdot \bmf(\bmx_I)$ and $h_I = \hat{\bmn} (\bmx_I) \cdot \bmf (\bmx_I) - \nu \hat{\bmn} (\bmx_I) \cdot \left( \nabla \times \nabla \times \bmv (\bmx_I) \right)$, using the form of the coefficient $\bmc$ from equation (\ref{eq:neumann-solution-coefficient}), the Laplacian operator for pressure field along the boundary $\Gamma$ will be approximated as:
\begin{align}
    \left( \nabla^2 \bmp (\bmx) |_{\bmx = \bmx_I} \right) \bmc (\bmx_I) &= s_I, \quad \text{ if } I \in I_\Gamma, \\
    \left( \nabla^2 \bmp (\bmx) |_{\bmx = \bmx_I}\right) \left( \bmR_{\bmvarphi \bmvarphi}(\bmx_I) \bmvarphi + \bmR_{\bmvarphi h}(\bmx_I) h_I \right) &= s_I, \quad \text{ if } I \in I_\Gamma.
\end{align}
This rearrangement allows the forcing term to remain on the right:
\begin{equation}
    \left( \nabla^2 \bmp(\bmx) |_{\bmx = \bmx_I}\right) \bmR_{\bmvarphi \bmvarphi}(\bmx_I) \bmvarphi = s_I - \left( \nabla^2 u(\bmx) |_{\bmx = \bmx_I}\right) \bmR_{\bmvarphi h}(\bmx_I) h_I, \quad \text{ if } I \in I_\Gamma, \label{eq:stokes-pressure-neumann-2}
\end{equation}
so that the left hand side only consists of coefficients of $p$. Now let $\beta_I = \left( \nabla^2 \bmp(\bmx) |_{\bmx = \bmx_I}\right) \bmR_{\bmvarphi h}(\bmx_I)$, equation (\ref{eq:stokes-pressure-neumann-2} can be re-written into:
\begin{equation}
    \sum_{\text{supp}(W_{Ij})} \bmalpha^2_{Ij} \varphi_j = \nabla \cdot \bmf(\bmx_I) - \beta_I h_I , \quad \text{ for } I \in I_\Gamma, \label{eq:grad-pressure-stencil-1}
\end{equation}
or, using the stencil for the velocity field in equation (\ref{eq:curl-vector-stencil}):
\begin{equation}
    \beta_I \hat{\bmn}(\bmx_I) \cdot \sum_{j \in \text{supp}(W_{Ij})} \bmalpha^3_{Ij} \bmv_j + \sum_{j \in \text{supp}(W_{Ij})} \bmalpha^1_{Ij} \varphi_j = \nabla \cdot \bmf(\bmx_I) - \beta_I \hat{\bmn} (\bmx_I) \cdot \bmf (\bmx_I) \quad \text{ for } I \in I_\Gamma.
\end{equation}

The following 2-by-2 block matrix system can then be assembled:
\begin{equation}
    \left[
        \begin{array}{c;{2pt/2pt}c;{2pt/2pt}c}
            \bmK & \bmG & \bm{0}^T \\ \hdashline[2pt/2pt]
            \bmN & \bmL & \bm{1}^T \\ \hdashline[2pt/2pt]
            \bm{0} & \bm{1} & N
        \end{array}
    \right]
    \left[
        \begin{array}{c}
            \bmv \\ \hdashline[2pt/2pt]
            \bmvarphi \\ \hdashline[2pt/2pt]
            \lambda
         \end{array}
    \right]
    =
    \left[
        \begin{array}{c}
            \bmb \\ \hdashline[2pt/2pt]
            \bmg \\ \hdashline[2pt/2pt]
            0
         \end{array}
    \right]. \label{eq:stokes-block-system}
\end{equation}
Each block is given as follows:
\begin{alignat}{2}
    \bmK_{IJ} &=
    \begin{cases}
        \delta_{IJ}, \quad &\text{ for } I \in I_\Gamma, \\
        \bmalpha^3_{Ij}, \quad &\text{ for } I \in I_\Omega, \\
    \end{cases} \quad \quad
    \bmG_{IJ} &&=
    \begin{cases}
        0, \quad &\text{ for } I \in I_\Gamma, \\
        \bmalpha^1_{Ij}, \quad &\text{ for } I \in I_\Omega, \\
    \end{cases} \quad \quad
    \bmN_{IJ} =
    \begin{cases}
        \hat{\bmn}_I \cdot \bmalpha^3_{Ij}, \quad &\text{ for } I \in I_\Gamma, \\
        0, \quad &\text{ for } I \in I_\Omega, \\
    \end{cases} \\
    \bmb_{I} &=
    \begin{cases}
        \bmw(\bmx_I), \quad &\text{ for } I \in I_\Gamma, \\
        \bmf(\bmx_I), \quad &\text{ for } I \in I_\Omega, \\
    \end{cases} \quad \quad
    \bmg_{I} &&=
    \begin{cases}
        \nabla \cdot \bmf(\bmx_I) - \beta_I \hat{\bmn}_I \cdot \bmf(\bmx_I), \quad &\text{ for } I \in I_\Gamma, \\
        \nabla \cdot \bmf(\bmx_I), \quad &\text{ for } I \in I_\Omega, \\
    \end{cases}
\end{alignat}
where a local-to-global index mapping is used. To be specific, the following equation:
\begin{equation}
    J = \text{LTG}(I, j),
\end{equation}
represents the mapping from the $j^{\text{th}}$ local index of particle $I$ to its respective global index $J$. Since the pressure field is effectively the Poisson problem with pure Neumann boundary condition, the solution contains a constant vector that must be accounted for in order to obtain a unique solution. Here, a zero-mean pressure field is enforced by adding a single Lagrange multiplier to the Poisson problem. This gives rise to the last row and column in equation (\ref{eq:stokes-block-system}). The $\bm{0}$ and $\bm{1}$ blocks are $N$-by-1 row vectors of $0$s and $1$s, respectively.

\subsection{Implementation}

The solution is implemented in \texttt{COMPADRE} \citep{paul_kuberry_2020_3876465}. There are two levels of parallelism in the implementation. The top level consists of domain decomposition through \texttt{ZOLTAN} package \citep{zoltan}, which splits the simulation domain into smaller domains. These domains are allocated to computing nodes and will be communicated through \texttt{MPI}.

At each node, parallel threads are executed to solve the least-squares problems at the local level. The local GMLS problem is solved using \texttt{LAPACK} package \citep{lapack}. The parallel execution of code and data management is implemented using the \texttt{KOKKOS} library \citep{edwards2014kokkos}. For the block matrix system, block Gauss-Seidel is used as the preconditioner through the \texttt{MueLu} package \citep{prokopenko2014muelu} inside the \texttt{Trilinos} project \citep{heroux2005overview}. It is noted that this is not the optimal choice, and will be addressed further in the discussion.

\section{Results}\label{results}

The simulation domain used for these tests is a bi-unit cube with dimension $[-1, 1] \times [-1, 1] \times [-1, 1]$. The particles are placed uniformly and equally distant from each other across the cube, and boundary particles lie on the six outer faces. Let $N$ be the number of particles along each direction. Therefore, there are $N \times N \times N$ particles across the domain. For the following results, the dynamic viscosity is set to be $1$.

\subsection{Manufactured solution: polynomial case}
The first cases will focus on manufacturing a solution that can be exactly reproduced from the basis functions. It is expected that the solution will be computed within numerical precision error. For this test, the manufactured solution when solving with quadratic basis functions is:
\begin{equation}
    \bmv = \begin{pmatrix}
        7x^2 + 6y^2 \\
        -6yz \\
        3z^2 - 14xz
    \end{pmatrix}, \quad \quad
    \varphi = x(1 + x + y) + y(1 + y + z) + z(1 + z + z^2), \\
    \label{eq:manufactured-solution-2nd}
\end{equation}
and for the quartic basis functions is:
\begin{equation}
    \bmv = \begin{pmatrix}
        7xz^2 + 6y^2 \\
        -7yz^2 \\
        -2x^3
    \end{pmatrix}, \quad \quad
    \varphi = x(x^2 + y^2 + z^2 + xy + yz + xz) + y(y^2 + z^2 + yz) + z^3. \\
    \label{eq:manufactured-solution-4th}
\end{equation}
In both cases, the velocity is chosen to be divergence-free. Using this implementation, the root mean square errors are collected in table \ref{tab:quadratic-polynomial} and table \ref{tab:quartic-polynomial}, respectively.
\begin{table}
    \parbox{0.5\linewidth}{
        \centering
        \begin{tabular}{|l|c|c|}
            \hline
               & \multicolumn{2}{c|}{$L_2$ error norm} \\ \hline
               N  & Velocity       & Pressure       \\ \hline
               6  & 6.80E-12    & 2.35E-11    \\ \hline
               12 & 2.84E-11    & 1.28E-10    \\ \hline
               24 & 2.00E-11    & 2.10E-10    \\ \hline
               48 & 4.90E-11    & 1.25E-09    \\ \hline
               96 & 8.51E-11    & 2.12E-09    \\ \hline
        \end{tabular}
        \caption{Quadratic polynomials.}
        \label{tab:quadratic-polynomial}
    } \hfill
    \parbox{0.5\linewidth}{
        \centering
        \begin{tabular}{|l|c|c|}
            \hline
               & \multicolumn{2}{c|}{$L_2$ error norm} \\ \hline
               N  & Velocity       & Pressure       \\ \hline
               6  & 1.60E-12    & 2.03E-10    \\ \hline
               12 & 3.50E-11    & 6.40E-10    \\ \hline
               24 & 5.90E-11    & 4.62E-10    \\ \hline
               48 & 7.29E-11    & 9.22E-09    \\ \hline
               96 & 6.18E-11    & 7.64E-09    \\ \hline
        \end{tabular}
        \caption{Quartic polynomial.}
        \label{tab:quartic-polynomial}
    }
\end{table}

\subsection{Manufactured solution: trigonometric case}

\begin{figure}[ht]
  \centering
  \includegraphics[width=\textwidth]{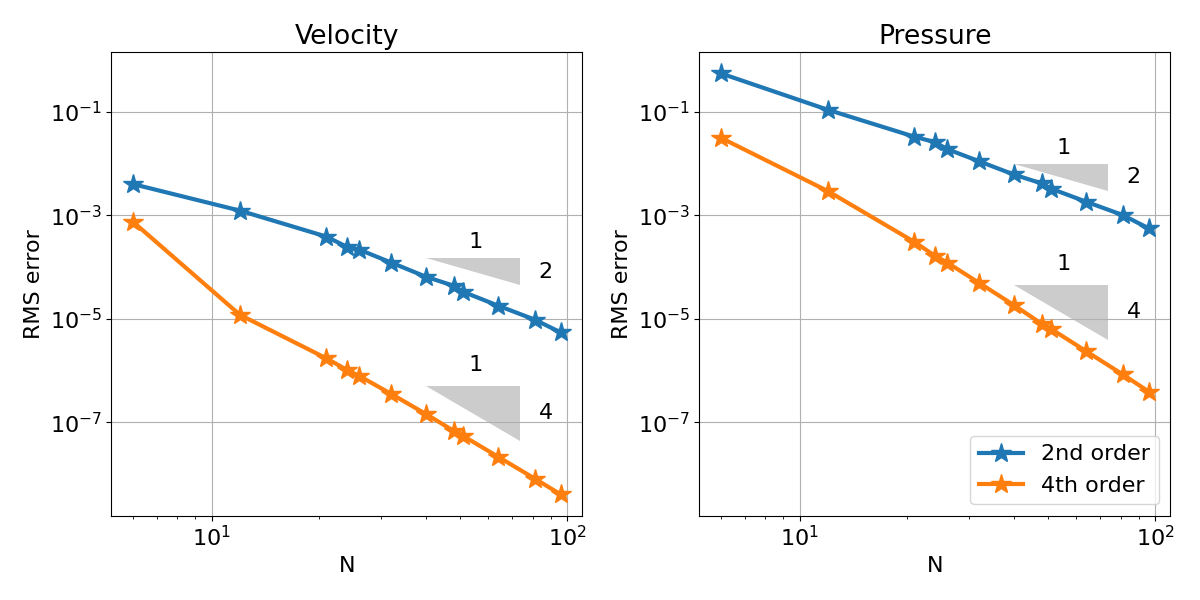}
  \caption{Convergence study for manufactured Stokes solution in 3D.}
    \label{fig:stokes-convergence-3D}
\end{figure}

In this test, the manufactured solution for velocity and pressure is as follows:
\begin{equation}
    \bmv = \begin{pmatrix}
        \sin{y} \sin {z} \\
        \sin{x} \sin {z} \\
        \sin{x} \sin {y} \\
    \end{pmatrix}, \quad \quad p = \sin{x} \sin{y} \sin{z}
    \label{eq:manufactured-solution-sin}
\end{equation}

Here, the reconstructed field does not lie in the space spanned by the basis functions. Therefore, it is expected that the $L_2$ error for the solution will converge with the same order of the basis functions used. More specifically, solutions obtained from using quadratic basis functions will show second-order error convergence, while quartic basis functions will show fourth order convergence for the root mean square error. As more particles are added to the simulation, and thus reducing the particles' spacing, the errors shown in Figure \ref{fig:stokes-convergence-3D} converge as expected in the asymptotic regime.

\subsection{Weak scaling performance}

\begin{figure}[ht]
  \centering
  \includegraphics[width=\textwidth]{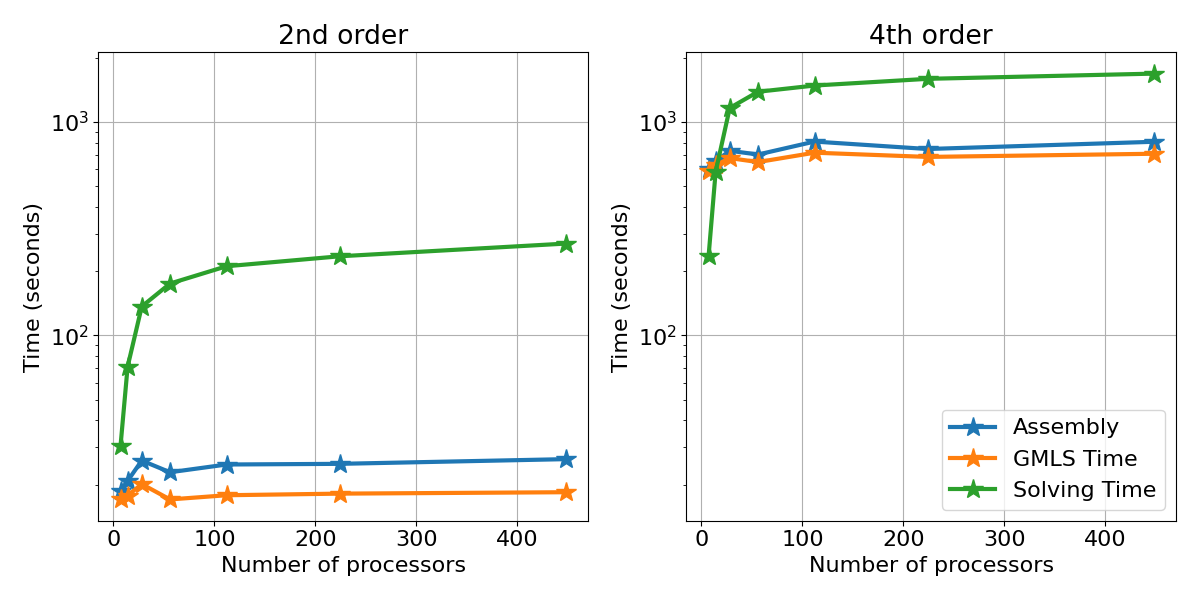}
  \caption{Weak scaling performance for the manufactured solution from equation (\ref{eq:manufactured-solution-sin}) with \texttt{COMPADRE}, using approximately 1200 particles per processor.}
  \label{fig:stokes-weak-scaling}
\end{figure}

In this section, the potential scalability of this approach is demonstrated. The ratio between degrees of freedom as number of processors used is kept constant as more processors are deployed to solve the Stokes problem. As seen in Figure \ref{fig:stokes-weak-scaling}, the time taken for assembling the block matrix system and solving least-squares GMLS problems remain quite constant in the weak scaling study. This is because these tasks are executed at the thread level, which is well performed when implemented with \texttt{KOKKOS} packages. The solving time, however, does not follow the same trend. The solve time increases slowly with problem size despite the concomitant increasing number of cores. This can be partially attributed to the sub-optimal choice of preconditioner.

\subsection{Trade-off between solving time and accuracy}

\begin{figure}[ht]
  \centering
  \includegraphics[width=\textwidth]{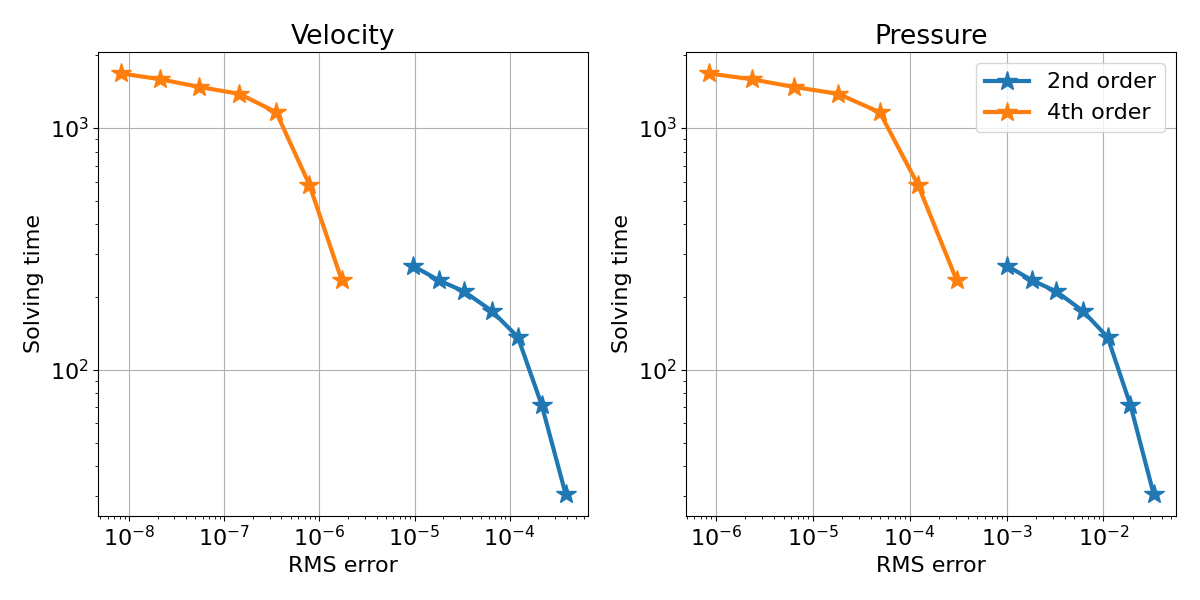}
  \caption{Trade-off between solving time and RMS error for second and fourth order basis. For a desired error below $10^{-5}$, it is more efficient to use a fourth order basis}
  \label{fig:stokes-timing}
\end{figure}

Figure \ref{fig:stokes-timing} shows the trade-off between solving time and accuracy between using quadratic and quartic polynomial basis. Because quartic polynomial basis consists of larger local least-squares problems, there are more non-zero entries in the global block system of the assembled Stokes solver, especially in the off-diagonal ones. Consequently, it takes longer to solve the Stokes flow with GMLS using quartic polynomial basis functions than quadratic ones. However, as the resolution enters the asymptotic regime in figure \ref{fig:stokes-timing}, the slope for the 4th order solver is significantly smaller when compared to the one obtained from the 2nd order solver. This suggests that the trade-off gain between solving time and numbers of significant figures is better for GMLS with quartic polynomials.

\subsection{Analytic solution - flow around sphere}

\begin{figure}[ht]
    \centering
    \includegraphics[width=0.6\textwidth]{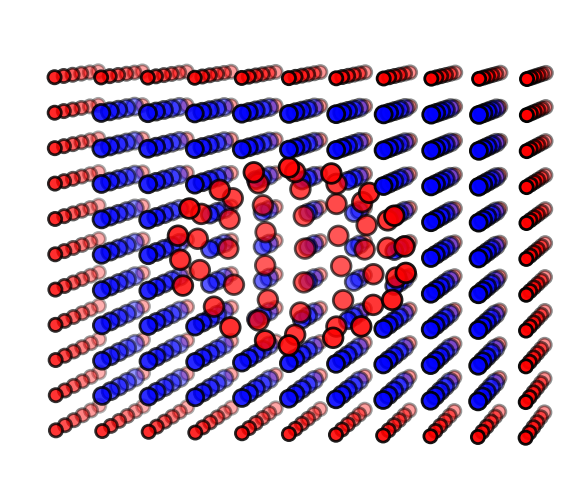}
    \caption{Sliced view of the simulation domain for Stokes' test case.}
    \label{fig:simulation-domain}
\end{figure}

One of the fundamental results in low Reynolds hydrodynamics is the Stokes' solution for steady flow past a small sphere. Let $a=1$ be the radius of the sphere placed at the center $(0, 0, 0)$ inside the box $[-2,2] \times [-2, 2] \times [-2, 2]$. Particles on the boundary of the sphere and the surfaces of the box are labeled red in Figure \ref{fig:simulation-domain}. The internal particles, colored blue in Figure \ref{fig:simulation-domain}, are placed uniformly in the space between the boundary surfaces. These particles are placed with average distance $h$ away from each other. In order to define the analytical solution, a spherical coordinate system $(r, \theta, \phi)$ is used, in which the origin $(0, 0, 0)$ is placed at the sphere's center. Let $W$ be the magnitude of the upward ambient velocity along the polar axis, the boundary condition for the velocity field $\bmv = (v_r, v_\theta, v_\phi)$ is:
\begin{alignat}{2}
    v_r &= v_\phi = 0, \quad &&\text{ at } r = a, \\
    v_r &= W \cos{\theta},  v_\theta = -W\sin{\theta}, \quad &&\text{ at } r = \infty.
\end{alignat}
With that, the analytical solution for the components in the fluid is:
\begin{align}
    v_r &= W \cos{\theta} \left(1 + \frac{a^3}{2r^3} - \frac{3a}{2r} \right), \\
    v_\theta &= -W \sin{\theta} \left(1 - \frac{a^3}{4r^3} - \frac{3a}{4r} \right), \\
    v_\phi &= 0,
\end{align}
and the gradient of the pressure field is:
\begin{equation}
    \nabla p = \frac{3}{2} \frac{\nu W a}{r^3} \cos{\theta}.
\end{equation}
In the following results, $\nu = 1$, $W = 10$ and $a =1$.

\begin{figure}[ht]
    \centering
    \includegraphics[width=\textwidth]{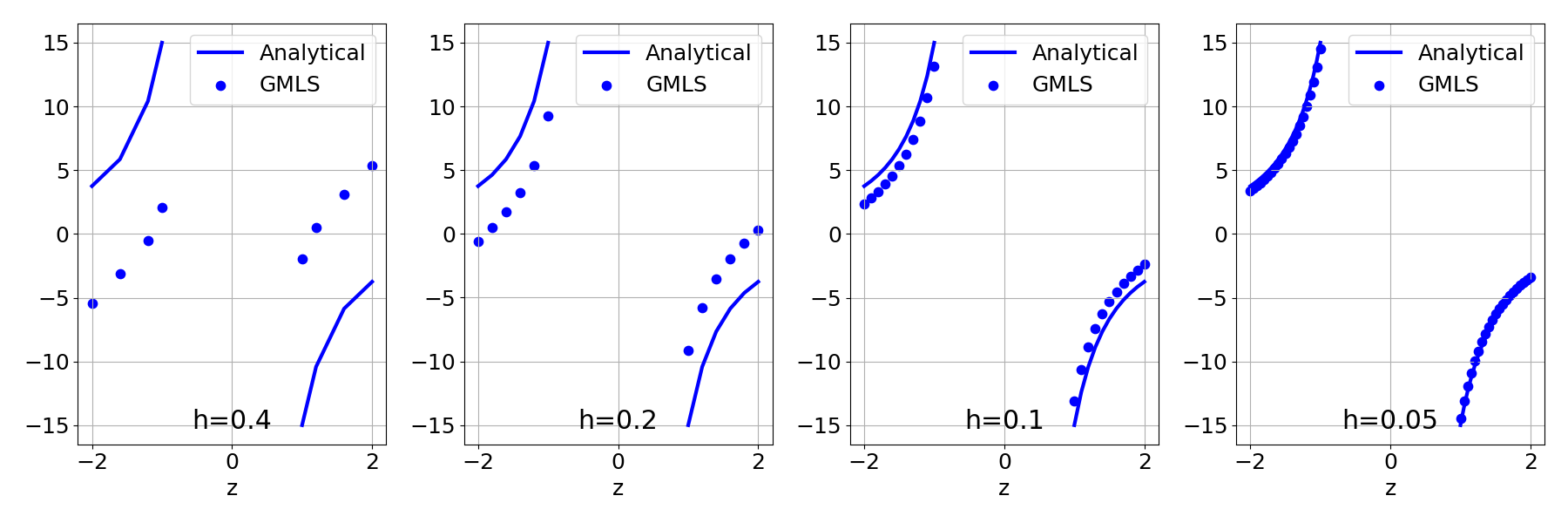}
    \caption{Plot of pressure across the line $x=0, y=0$.}
    \label{fig:stokes-sphere-pressure-line}
\end{figure}

\begin{figure}[ht]
    \centering
    \includegraphics[width=\textwidth]{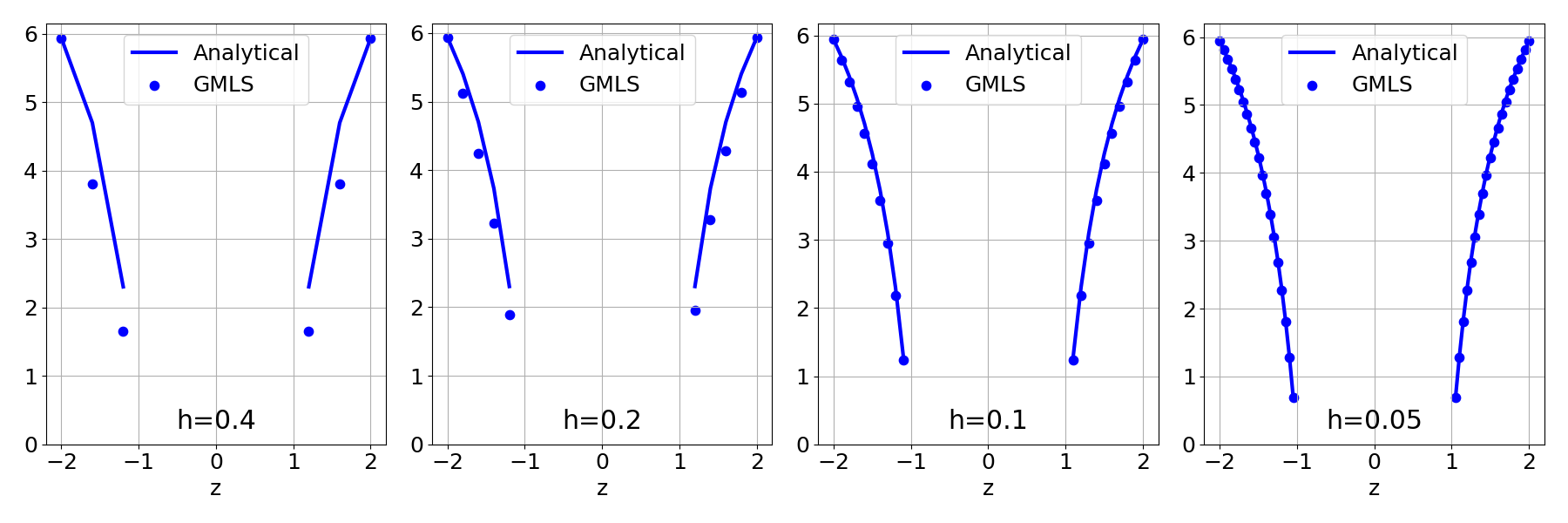}
    \caption{Plot of velocity's magnitude across the line $x=0, y=0$.}
    \label{fig:stokes-sphere-velocity-line}
\end{figure}

\begin{figure}[ht]
    \centering
    \includegraphics[width=\textwidth]{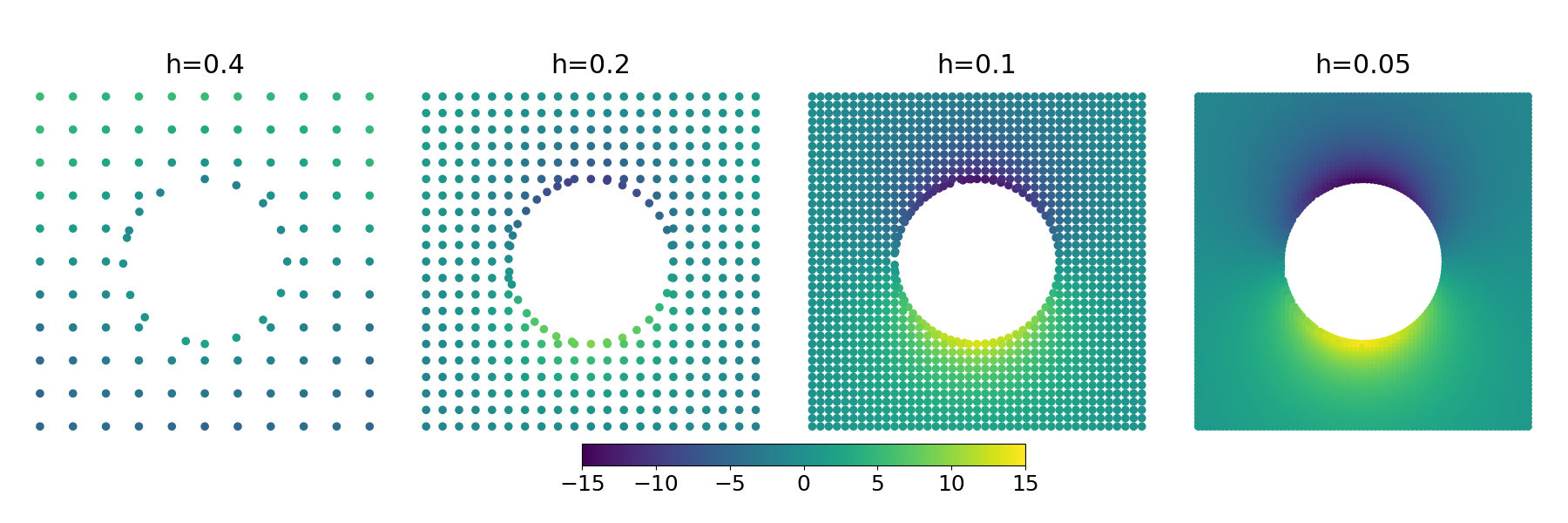}
    \caption{Scatter plot of pressure in the plane $y=0$.}
    \label{fig:stokes-sphere-pressure-plane}
\end{figure}

\begin{figure}[ht]
    \centering
    \includegraphics[width=\textwidth]{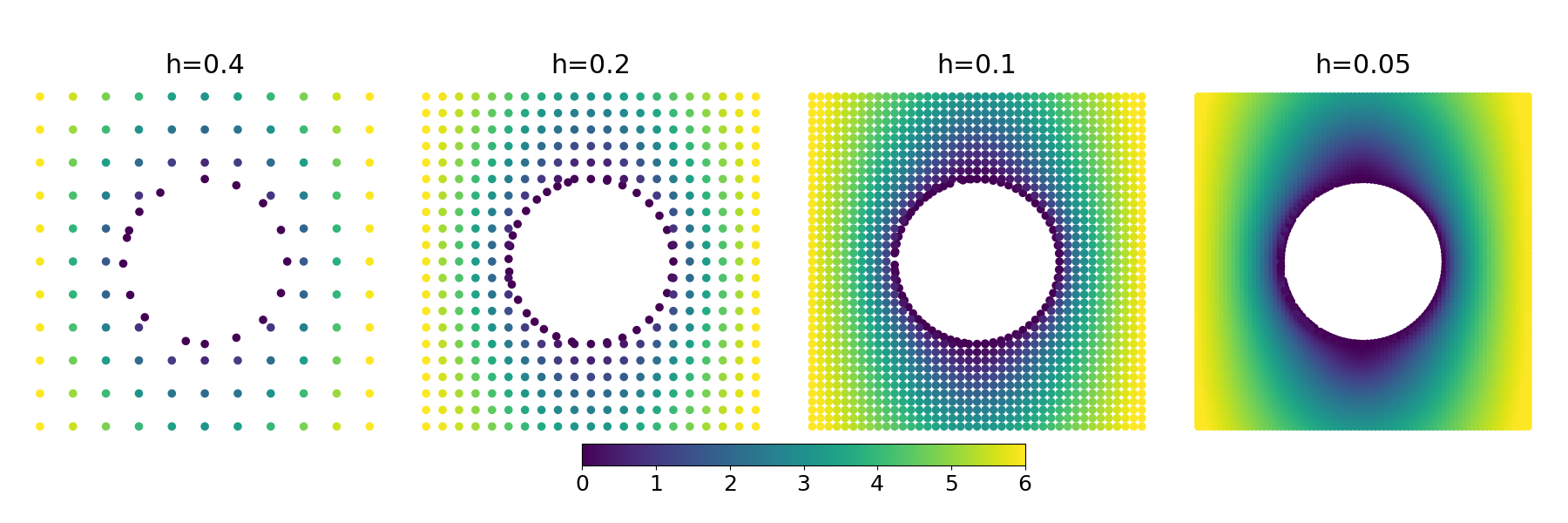}
    \caption{Scatter plot of velocity's magnitude in the plane $y=0$.}
    \label{fig:stokes-sphere-velocity-plane}
\end{figure}

\begin{figure}[ht]
  \centering
  \includegraphics[width=\textwidth]{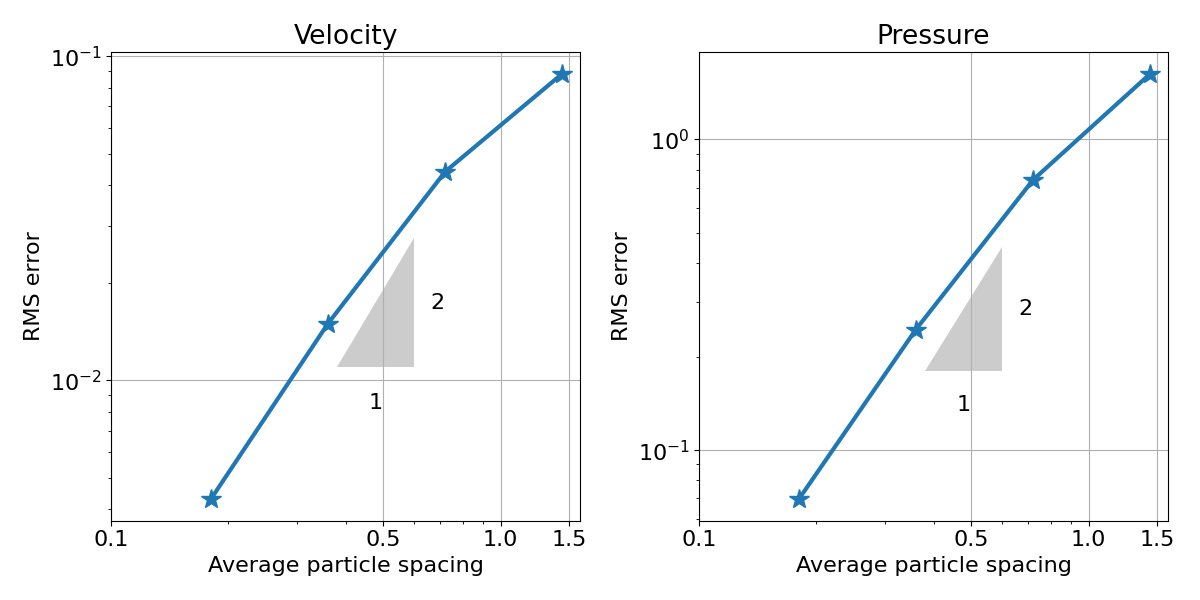}
  \caption{Convergence rate for refinement study of Stokes' flow around sphere with quadratic GMLS.}
    \label{fig:stokes-sphere-convergence}
\end{figure}

The computed result is compared to the analytical formula. Figures \ref{fig:stokes-sphere-pressure-line} and \ref{fig:stokes-sphere-velocity-line} compare the values of the computed and analytical fields along the axis line of $x = y =0$. It is apparent that the computation provides better results as the simulation domain gets refined (the particle's spacing $h$ decreases). This is also observed from the contour plot on the plane $y=0$, as illustrated on figures \ref{fig:stokes-sphere-pressure-plane} and \ref{fig:stokes-sphere-velocity-plane}. The convergence rate for the root mean square error for this case is shown in figure \ref{fig:stokes-sphere-convergence}.

\section{Conclusion}\label{conclusion}
This work has presented the software implementation solving the steady state Stokes flow with generalized moving least squares in \texttt{COMPADRE}. A previous paper \citep{Trask2018} presented for the first time a fully meshless method for the Stokes problem that is able to achieve high-order convergence for both the velocity and the pressure, while maintaining a sparse discretization. The current paper shows that the behaviour demonstrated in \citep{Trask2018} extends to 3D, and that the methodology lends itself to scalable parallel implementation with encouraging scaling behaviour for large problems. Unlike other mesh-based methods, this meshless approach does not require computationally expensive mesh generation. The cost of generating a compatible point cloud is negligible when compared to the cost of solving the entire system of partial differential equations. After having showed the accuracy and efficiency of the approach, analytical solutions are used to systematically benchmark the necessary components to use the scheme to study steady-state Stokes flow. 

The great performance demonstrated in the weak scaling study shows the promising scalability for this implementation, as the numerical solver has two parallelism levels. The domain decomposition, achieved through \texttt{Zoltan2}, is allocated to each computing node and communication between them is handled through \texttt{MPI}. The local least-square problems are allocated to threads in each node, which use the shared memory on their respective parent node. The parallelism at the thread level is handled through \texttt{KOKKOS}. Block Gauss-Seidel preconditioner is applied through \texttt{MueLu} package inside the \texttt{Trilinos} project. It is also worth noting that the preconditioner used for solving the block system is far from the optimal one. It is expected that the solving time would improve once an optimal choice is used, and further discussion on this matter will be saved for future work.

\section*{Acknowledgments}

Q.-T. Ha and E.M. Ryan acknowledge the support of the National Science Foundation through grants 1727316 and 1911698.

Sandia  National  Laboratories  is  a  multi-mission  laboratory  managed  and  operated  by  National  Technology and Engineering Solutions of Sandia, LLC., a wholly owned subsidiary of Honeywell International,Inc.,  for  the  U.S.  Department  of  Energy’s  National  Nuclear  Security  Administration  under  contract  DE-NA0003525.  This paper describes objective technical results and analysis.  Any subjective views or opinions that might be expressed in the paper do not necessarily represent the views of the U.S. Department of Energy or the United States Government. 

N. Trask’s work is supported under the U.S. Department of Energy, Office of Science, Office of Advanced Scientific Computing  Research  under  the  Collaboratory  on  Mathematics  and  Physics-Informed  Learning Machines for Multiscale and Multiphysics Problems (PhILMs) project.

P. Kuberry's work is supported by the Coupling Approaches for Next-Generation Architectures (CANGA) project, a joint effort funded by the United States Department of Energy’s Office of Science under the Scientific Discovery through Advanced Computing (SciDAC), and the Laboratory Directed Research and Development program at Sandia National Laboratories, United States.

SAND number: SAND2021-0727 O


\subsection*{Financial disclosure}

None reported.

\subsection*{Conflict of interest}

The authors declare no potential conflict of interests.

\bibliographystyle{plainnat}
\bibliography{references}  






\end{document}